\documentclass[10pt,twoside]{article}
\usepackage{graphicx}
\usepackage{amsmath}
\usepackage{amssymb}
\usepackage{amscd}
\usepackage{Latex-document}

\newenvironment{enumerateA}{%
 \begin{enumerate}}{\end{enumerate}}

\DeclareMathOperator{\car}{char} \DeclareMathOperator{\Ch}{CH}
\DeclareMathOperator{\Fix}{Fix} \DeclareMathOperator{\Gal}{Gal}
\DeclareMathOperator{\Nrd}{Nrd} \DeclareMathOperator{\SL}{SL}
\DeclareMathOperator{\SO}{SO} \DeclareMathOperator{\Spec}{Spec}

\newcommand{\amod}{{\rm mod}\mkern5mu} 
\newcommand{\et}{\mathrm{\acute et}}
\newcommand{\degre}[2]{[#1\mathbin{\mathop:}#2]}

\newcommand{\LC}{\mathbf C}
\newcommand{\LN}{\mathbf N}
\newcommand{\LP}{\mathbf P}
\newcommand{\LZ}{\mathbf Z}

\newcommand{\bmidb}[2]{\{\,#1\mid#2\,\}}
\newcommand{\Bmidb}[2]{\left\{\,#1
    \mathrel{\left.\vphantom{#1}\vphantom{#2}\right|}#2\,\right\}}
\newcommand{\amida}[2]{\langle\,#1\mid#2\,\rangle}

\newcommand{\pfister}[1]{\langle\!\langle#1\rangle\!\rangle}
\newcommand{\qform}[1]{\langle#1\rangle}

\title{\bf Norm Varieties and Algebraic Cobordism\vskip 6mm}
\author{Markus Rost\thanks{Department of Mathematics, The Ohio State University, 231 W 18th Avenue,
Columbus, OH 43210, USA. E-mail: rost@math.ohio-state.edu, URL: http://www.math.ohio-state.edu/\~{}rost}}

\date{\vspace{-8mm}}

\begin{document}

\maketitle

\thispagestyle{first} \setcounter{page}{77}

\begin{abstract}\vskip 3mm

We outline briefly results and examples related with the bijectivity of the norm residue homomorphism.  We define
norm varieties and describe some constructions.  We discuss degree formulas which form a major tool to handle norm
varieties.  Finally we formulate Hilbert's 90 for symbols which is the hard part of the bijectivity of the norm
residue homomorphism, modulo a theorem of Voevodsky.

\vskip 4.5mm

\noindent {\bf 2000 Mathematics Subject Classification:} 12G05.

\noindent {\bf Keywords and Phrases:} Milnor's $K$-ring, Galois
cohomology, Cobordism.
\end{abstract}

\vskip 12mm

\section*{Introduction}

\vskip-5mm \hspace{5mm}

This text is a brief outline of results and examples related with
the bijectivity of the norm residue homomorphism---also called
``Bloch-Kato conjecture'' and, for the $\amod 2$ case, ``Milnor
conjecture''.

The starting point was a result of Voevodsky which he communicated
in 1996.  Voevodsky's theorem basically reduces the Bloch-Kato
conjecture to the existence of norm varieties and to what I call
Hilbert's 90 for symbols.  Unfortunately there is no text
available on Voevodsky's theorem.

In this exposition $p$ is a prime, $k$ is a field with $\car k\neq
p$ and $K^M_nk$ denotes Milnor's $n$-th $K$-group of~$k$
\cite{Milnor:70}, \cite{Rost:96}.

Elements in $K^M_nk/p$ of the form
\begin{displaymath}
  u=\{a_1,\dots,a_n\}\bmod p
\end{displaymath}
are called symbols ($\amod p$, of weight~$n$).

A field extension $F$ of~$k$ is called a splitting field of $u$ if
$u_F=0$ in $K^M_nF/p$.

Let
\begin{gather*}
  h_{(n,p)}\colon K^M_nk/p\to H^n_\et(k,\mu_p^{\otimes n}), \\
  \{a_1,\dots,a_n\}\mapsto (a_1,\dots,a_n)
\end{gather*}
be the norm residue homomorphism.

\section{Norm varieties}

\vskip-5mm \hspace{5mm}

All successful approaches to the Bloch-Kato conjecture consist of
an investigation of appropriate generic splitting varieties of
symbols. This goes back to the work of Merkurjev and Suslin on the
case $n=2$ who studied the $K$-cohomology of Severi-Brauer
varieties \cite{Merkurjev-Suslin:82}.  Similarly, for the case
$p=2$ (for $n=3$ by Merkurjev, Suslin \cite{Merkurjev-Suslin:86}
and the author \cite{Rost:1986x}, for all $n$ by Voevodsky
\cite{Voevodsky:02}) one considers certain quadrics associated
with Pfister forms.  For a long time it was not clear which sort
of varieties one should consider for arbitrary $n$,~$p$.  In some
cases one knew candidates, but these were non-smooth varieties and
desingularizations appeared to be difficult to handle.  Finally
Voevodsky proposed a surprising characterization of the necessary
varieties.  It involves characteristic numbers and yields a
beautiful relation between symbols and cobordism theory.

{\bf Definition.} \it
  Let $u=\{a_1,\dots,a_n\}\bmod p$ be a symbol.  Assume that $u\neq
  0$.  A \emph{norm variety} for~$u$ is a smooth proper irreducible
  variety~$X$ over~$k$ such that
  \begin{enumerateA}
    \item
    \label{item:split}
    The function field~$k(X)$ of~$X$ splits $u$.
    \item
    $\dim X=d:=p^{n-1}-1$.
    \item
    \label{item:sd}
    $\frac{s_d(X)}p\not\equiv 0\mod p$.
  \end{enumerateA}
\rm

Here $s_d(X)\in\LZ$ denotes the characteristic number of~$X$ given
by the $d$-th Newton polynomial in the Chern classes of~$TX$.  It
is known (by Milnor) that in dimensions $d=p^n-1$ the
number~$s_d(X)$ is $p$-divisible for any~$X$.  If $k\subset \LC$
one may rephrase condition~\eqref{item:sd} by saying that $X(\LC)$
is indecomposable in the complex cobordism ring $\amod p$.

We will observe in section~\ref{sec:degree-formulas} that the
conditions for a norm variety are birational invariant.

The name ``norm variety'' originates from some constructions of
norm varieties, see section~\ref{sec:exist-norm-vari}.

We conclude this section with the ``classical'' examples of norm
varieties.

{\bf Example.}
  The case $n=2$.  Assume that $k$ contains a primitive $p$-th
  root~$\zeta$ of unity.  For $a$, $b\in k^*$ let $A_\zeta(a,b)$ be
  the central simple $k$-algebra with presentation
  \begin{displaymath}
    A_\zeta(a,b)=\amida{u,v}{u^p=a,v^p=b,vu=\zeta uv}.
  \end{displaymath}
  The Severi-Brauer variety $X(a,b)$ of~$A_\zeta(a,b)$ is a norm
  variety for the symbol $\{a,b\}\bmod p$.

{\bf Example.}
  The case $p=2$.  For $a_1$, \ldots, $a_n\in k^*$ one denotes by
  \begin{displaymath}
    \pfister{a_1,\ldots,a_n}=\bigotimes_1^n\qform{1,-a_i},
  \end{displaymath}
  the associated $n$-fold Pfister form \cite{Lam:73},
  \cite{Scharlau:85}.  The quadratic form
  \begin{displaymath}
    \varphi=\pfister{a_1,\ldots,a_{n-1}}\perp\qform {-a_n}
  \end{displaymath}
  is called a Pfister neighbor.  The projective quadric $Q(\varphi)$
  defined by $\varphi=0$ is a norm variety for the symbol
  $\{a_1,\ldots,a_{n}\}\bmod 2$.

\section{Degree formulas} \label{sec:degree-formulas}

\vskip-5mm \hspace{5mm}

The theme of ``degree formulas'' goes
back to Voevodsky's first text on the Milnor conjecture (although
he never formulated explicitly a ``formula'') \cite{Voevodsky:96}.
In this section we formulate the degree formula for the
characteristic numbers~$s_d$.  It shows the birational invariance
of the notion of norm varieties.

The first proof of this formula relied on Voevodsky's stable
homotopy theory of algebraic varieties.  Later we found a rather
elementary approach \cite{Merkurjev:2000x}, which is in spirit
very close to ``elementary'' approaches to the complex cobordism
ring \cite{Quillen:71}, \cite{Buoncristiano-Hacon:83}.

For our approach to Hilbert's 90 for symbols we use also ``higher
degree formulas'' which again were first settled using Voevodsky's
stable homotopy theory~\cite{Borghesi:00}.  These follow meanwhile
also from the ``general degree formula'' proved by Morel and
Levine \cite{Levine-Morel:02} in characteristic~$0$ using
factorization theorems for birational maps \cite{MR1896232}.

We fix a prime $p$ and a number~$d$ of the form $d=p^{n}-1$.

For a proper variety~$X$ over~$k$ let
\begin{displaymath}
  I(X)=\deg\bigl(\Ch_0(X)\bigr)\subset\LZ
\end{displaymath}
be the image of the degree map on the group of $0$-cycles.  One
has $I(X)=i(X)\LZ$ where $i(X)$ is the ``index'' of~$X$, i.~e.,
the gcd of the degrees $\degre{k(x)}k$ of the residue class field
extensions of the closed points~$x$ of~$X$.  If $X$ has a
$k$-point (in particular if $k$ is algebraically closed), then
$I(X)=\LZ$.  The group~$I(X)$ is a birational invariant of~$X$.
We put
\begin{displaymath}
  J(X)=I(X)+p\LZ.
\end{displaymath}

Let $X$, $Y$ be irreducible smooth proper varieties over~$k$ with
$\dim Y=\dim X=d$ and let $f\colon Y\to X$ be a morphism.  Define
$\deg f$ as follows: If $\dim f(Y)<\dim X$, then $\deg f=0$.
Otherwise $\deg f\in \LN$ is the degree of the extension
$k(Y)/k(X)$ of the function fields.

{\bf Theorem (Degree formula for \boldmath $s_d$).}
  \label{formula}
  \begin{displaymath}
    \frac{s_{d}(Y)}p=(\deg f)\frac{s_{d}(X)}p \mod J(X).
  \end{displaymath}

{\bf Corollary.} \it
  \label{birat}
  The class
  \begin{displaymath}
    \frac{s_{d}(X)}p\bmod J(X)\in \LZ/J(X)
  \end{displaymath}
  is a birational invariant.
\rm

{\bf Remark.}
  If $X$ has a $k$-rational point, then $J(X)=\LZ$ and the degree
  formula is empty.  The degree formula and the birational invariants
  $s_d(X)/p\bmod J(X)$ are phenomena which are interesting only over
  non-algebraically closed fields.  Over the complex numbers the only
  characteristic numbers which are birational invariant are the Todd
  numbers.

We apply the degree formula to norm varieties.  Let $u$ be a
nontrivial symbol $\amod p$ and let $X$ be a norm variety for~$u$.
Since $k(X)$ splits~$u$, so does any residue class field~$k(x)$
for $x\in X$.  As $u$ is of exponent~$p$, it follows that
$J(X)=p\LZ$.

{\bf Corollary (Voevodsky).} \it
  Let $u$ be a nontrivial symbol and let $X$ be a norm variety of~$u$.
  Let further $Y$ be a smooth proper irreducible variety with $\dim
  Y=\dim X$ and let $f\colon Y\to X$ be a morphism.  Then $Y$ is a
  norm variety for~$u$ if and only if $\deg f$ is prime to~$p$.
\rm

It follows in particular that the notion of norm variety is
birational invariant.  Therefore we may call any irreducible
variety~$U$ (not necessarily smooth or proper) a norm variety of a
symbol~$u$ if~$U$ is birational isomorphic to a smooth and proper
norm variety of~$u$.

\section{Existence of norm varieties}
\label{sec:exist-norm-vari}

\vskip-5mm \hspace{5mm}

{\bf Theorem.} \it
  Norm varieties exists for every symbol~$u\in K_n^Mk/p$ for every~$p$
  and every~$n$.
\rm

As we have noted, for the case $n=2$ one can take appropriate
Severi-Brauer varieties (if $k$ contains the $p$-th roots of
unity) and for the case $p=2$ one can take appropriate quadrics.

In this exposition we describe a proof for the case $n=3$ using
fix-point theorems of Conner and Floyd in order to compute the
non-triviality of the characteristic numbers.  Our first proof for
the general case used also Conner-Floyd fix-point theory.  Later
we found two further methods which are comparatively simpler.
However the Conner-Floyd fix-point theorem is still used in our
approach to Hilbert's~90 for symbols.

Let $u=\{a,b,c\}\bmod p$ with $a$, $b$,~$c\in k^*$.  Assume that
$k$ contains a primitive $p$-th root~$\zeta$ of unity, let
$A=A_\zeta(a,b)$ and let
\begin{displaymath}
  \mathit{MS}(A,c)=\bmidb{x\in A}{\Nrd(x)=c}.
\end{displaymath}
We call $\mathit{MS}(A,c)$ the Merkurjev-Suslin variety associated
with~$A$ and~$c$.  The symbol~$u$ is trivial if and only if
$\mathit{MS}(A,c)$ has a rational point
\cite{Merkurjev-Suslin:82}. The variety $\mathit{MS}(A,c)$ is a
twisted form of $\SL(p)$.

{\bf Theorem.} \it
  Suppose $u\neq0$.  Then $\mathit{MS}(A,c)$ is a norm variety
  for~$u$.
\rm

Let us indicate a proof for a subfield $k\subset \LC$ (and for
$p>2$). Let $U=\mathit{MS}(A,c)$.  It is easy to see that $k(U)$
splits~$u$. Moreover one has $\dim U=\dim A-1=p^2-1$.  It remains
to show that there exists a proper smooth completion~$X$ of~$U$
with nontrivial characteristic number.

Let
\begin{displaymath}
  \bar U=\bmidb{[x,t]\in \LP(A\oplus k)}{\Nrd(x)=ct^p}
\end{displaymath}
be the naive completion of~$U$.  We let the group $G=\LZ/p\times
\LZ/p$ act on the algebra~$A$ via
\begin{displaymath}
  (r,s)\cdot u=\zeta^r u, \quad (r,s)\cdot v=\zeta^s v.
\end{displaymath}
This action extends to an action on $\LP(A\oplus k)$ (with the trivial
action on~$k$) which induces a $G$-action on~$\bar U$. Let $\Fix(\bar
U)$ be the fixed point scheme of this action.  One finds that $\Fix(\bar
U)$ consists just of the $p$ isolated points $[1,\zeta^i \root p \of
c]$, $i=1$, \ldots,~$p$, which are all contained in~$U$.

The variety~$U$ is smooth, but~$\bar U$ is not.  However, by
equivariant resolution of singularities
\cite{Bierstone-Milman:97}, there exists a smooth proper
$G$-variety~$X$ together with a $G$-morphism $X\to \bar U$ which
is a birational isomorphism and an isomorphism over~$U$.  It
remains to show that
\begin{displaymath}
  \frac{s_d(X)}p\not\equiv 0\mod p.
\end{displaymath}
For this we may pass to topology and try to compute
$s_d\bigl(X(\LC)\bigr)$.  We note that for odd $p$, the Chern
number~$s_d$ is also a Pontryagin number and depends only on the
differentiable structure of the given variety.  Note further that
$X$ has the same $G$-fixed points as~$\bar U$ since the
desingularization took place only outside~$U$.

Consider the variety
\begin{displaymath}
  \textstyle Z=\Bmidb{\left[\sum_{i,j=1}^px_{ij}u^iv^j,t\right]\in
  \LP(A\oplus k)}{\sum_{i,j=1}^px_{ij}^p=c t^p}.
\end{displaymath}
This variety is a smooth hypersurface and it is easy to check
\begin{displaymath}
  \frac{s_d(Z)}p\not\equiv 0\mod p.
\end{displaymath}
As a $G$-variety, the variety~$Z$ has the same fixed points as~$X$
(``same'' means that the collections of fix-points together with
the $G$-structure on the tangent spaces are isomorphic).  Let $M$
be the differentiable manifold obtained from $X(\LC)$ and
$-Z(\LC)$ by a multi-fold connected sum along corresponding fixed
points.  Then $M$ is a $G$-manifold without fixed points.  By the
theory of Conner and Floyd \cite{Conner-Floyd:64}, \cite{Floyd:71}
applied to $(\LZ/p)^2$-manifolds of dimension $d=p^2-1$ one has
\begin{displaymath}
  \frac{s_d(M)}p\equiv 0\mod p.
\end{displaymath}
Thus
\begin{displaymath}
  \frac{s_d(X)}p\equiv\frac{s_d(Z)}p \mod p
\end{displaymath}
and the desired non-triviality is established.

{\bf The functions \boldmath$\Phi_n$.} We conclude this section with examples of norm varieties for the general
case.

Let $a_1$, $a_2$, \ldots be a sequence of elements in~$k^*$.  We
define functions $\Phi_n=\Phi_{a_1,\ldots,a_n}$ in~$p^n$ variables
inductively as follows.
\begin{align*}
  \Phi_0(t)&=t^p, \\
  \Phi_n(T_0,\ldots,T_{p-1})&=\Phi_{n-1}(T_0)\prod_{i=1}^{p-1}
  \bigl(1-a_n\Phi_{n-1}(T_i)\bigr).
\end{align*}
Here the~$T_i$ stand for tuples of $p^{n-1}$ variables.  Let
$U(a_1,\ldots,a_n)$ be the variety defined by
\begin{displaymath}
  \Phi_{a_1,\ldots,a_{n-1}}(T)=a_n.
\end{displaymath}

{\bf Theorem.} \it
  Suppose that the symbol $u=\{a_1,\dots,a_n\}\bmod p$ is nontrivial.
  Then $U(a_1,\ldots,a_n)$ is a norm variety of~$u$.
\rm

\section{Hilbert's 90 for symbols}
\label{sec:hilberts-90-symbols}

\vskip-5mm \hspace{5mm}

The bijectivity of the norm residue homomorphisms has always been
considered as a sort of higher version of the classical Hilbert's
Theorem~90 (which establishes the bijectivity for $n=1$).  In
fact, there are various variants of the Bloch-Kato conjecture
which are obvious generalizations of Hilbert's Theorem~90:  The
Hilbert's Theorem~90 for $K_n^M$ of cyclic extensions or the
vanishing of the motivic cohomology group
$H^{n+1}\bigl(k,\LZ(n)\bigr)$.  In this section we describe a
variant which on one hand is very elementary to formulate and on
the other hand is the really hard part of the Bloch-Kato
conjecture (modulo Voevodsky's theorem).

Let $u=\{a_1,\dots,a_n\}\in K^M_nk/p$ be a symbol.  Consider the
norm map
\begin{displaymath}
  \mathcal N_u=\sum_F N_{F/k}\colon \bigoplus_F K_1F\to K_1k
\end{displaymath}
where $F$ runs through the finite field extensions of~$k$
(contained in some algebraic closure of~$k$) which split~$u$.
Hilbert's Theorem 90 for~$u$ states that $\ker \mathcal N_u$ is
generated by the ``obvious'' elements.

To make this precise, we consider two types of basic relations
between the norm maps $N_{F/k}$.

Let $F_1$, $F_2$ be finite field extensions of~$k$.  Then the
sequence
\begin{equation}
  \label{eq:1}
  K_1(F_1\otimes F_2) \xrightarrow{ (N_{F_1\otimes F_2/F_1} ,
  -N_{F_1\otimes F_2/F_2}) } K_1F_1\oplus K_1F_2 \xrightarrow{
  N_{F_1/k} + N_{F_2/k}}K_1k
\end{equation}
is a complex.

Further, if $K/k$ is of transcendence degree~$1$, then the
sequence
\begin{equation}
  \label{eq:2}
  K_2 K\xrightarrow{d_K} \bigoplus_v K_1\kappa(v)\xrightarrow N K_1k
\end{equation}
is a complex.  Here $v$ runs through the valuations of $K/k$,
$d_K$ is given by the tame symbols at each~$v$ and $N$ is the sum
of the norm maps $N_{\kappa(v)/k}$.  The sum formula $N\circ
d_K=0$ is also known as Weil's formula.

We now restrict again to splitting fields of~$u$.  The maps
in~\eqref{eq:1} yield a map
\begin{displaymath}
  \mathcal R_u=\sum _{F_1,F_2} (N_{F_1\otimes F_2/F_1} ,
  -N_{F_1\otimes F_2/F_2})\colon \bigoplus_{F_1,F_2} K_1(F_1\otimes
  F_2) \to \bigoplus_F K_1F
\end{displaymath}
with $\mathcal N_u\circ \mathcal R_u=0$.  Let $C$ be the cokernel
of $\mathcal R_u$ and let $\mathcal N_u'\colon C\to K_1k$ be the
map induced by~$\mathcal N_u$.  Then the maps in \eqref{eq:2}
yield a map
\begin{displaymath}
  \mathcal S_u=\sum _{K} d_K\colon \bigoplus_K K_2K \to C
\end{displaymath}
with $\mathcal N_u'\circ \mathcal S_u=0$ where $K$ runs through
the splitting fields of~$u$ of transcendence degree~$1$ over~$k$
(contained in some universal field).  Let $ H_0(u,K_1)$ be the
cokernel of $\mathcal S_u$ and let $N_u\colon H_0(u,K_1)\to K_1k$
be the map induced by~$\mathcal N_u'$.

{\bf Hilbert's 90 for symbols.} \it
  For every symbol~$u$ the norm map
  \begin{displaymath}
    N_u\colon H_0(u,K_1)\to K_1k
  \end{displaymath}
  is injective.
\rm

{\bf Example.}
  If $u=0$, then it is easy to see that $N_u$ is injective.  In fact,
  it is a trivial exercise to check that $\mathcal N_u'$ is injective.

{\bf Example.}
  The case $n=1$.  The splitting fields~$F$ of $u=\{a\}\bmod p$ are
  exactly the field extensions of~$k$ containing a $p$-th root of~$a$.
  It is an easy exercise to reduce the injectivity of~$N_u$ (in fact
  of $\mathcal N_u'$) to the classical Hilbert's Theorem~90, i.~e.,
  the exactness of
  \begin{displaymath}
    K_1L\xrightarrow{1-\sigma}K_1L\xrightarrow{N_{L/k}}K_1k
  \end{displaymath}
  for a cyclic extension $L/k$ of degree~$p$ with~$\sigma$ a generator
  of $\Gal(L/k)$.

{\bf Example.}
  The case $n=2$.  Assume that $k$ contains a primitive $p$-th
  root~$\zeta$ of unity.  The splitting fields~$F$ of $u=\{a,b\}\bmod
  p$ are exactly the splitting fields of the algebra~$A_\zeta(a,b)$.
  One can show that
  \begin{displaymath}
    H_0(u,K_1)=K_1A_\zeta(a,b)
  \end{displaymath}
  with $N_u$ corresponding to the reduced norm map~$\Nrd$
  \cite{Merkurjev-Suslin:92}.  Hence in this case Hilbert's 90 for~$u$
  reduces to the classical fact $\mathit{SK}_1A=0$ for central simple
  algebras of prime degree \cite{Draxl:83}.

{\bf Example.}
  The case $p=2$.  The splitting fields~$F$ of
  $u=\{a_1,\dots,a_n\}\bmod 2$ are exactly the field extensions of~$k$
  which split the Pfister form $\pfister{a_1,\ldots,a_n}$ or,
  equivalently, over which the Pfister neighbor
  $\pfister{a_1,\ldots,a_{n-1}}\perp\qform {-a_n}$ becomes isotropic.
  Hilbert's 90 for symbols $\amod 2$ had been first established
  in~\cite{Ro:88}.  This text considered similar norm maps associated
  with any quadratic form (which are not injective in general).  A
  treatment of the special case of Pfister forms is contained in
  \cite{Kahn:97}.

{\bf Remark.}
  One can show that the group $H_0(u,K_1)$ as defined above is also
  the quotient of $\oplus_F K_1F$ by the $R$-trivial elements in $\ker
  \mathcal N_u$.  This is quite analogous to the description of $K_1A$
  of a central simple algebra~$A$:  The group $K_1A$ is the quotient
  of $A^*$ by the subgroup of $R$-trivial elements in the kernel of
  $\Nrd\colon A^*\to F^*$.  Similarly for the case $p=2$:  In this
  case the injectivity of~$N_u$ is related with the fact that for
  Pfister neighbors~$\varphi$ the kernel of the spinor norm
  $\SO(\varphi)\to k^*/(k^*)^2$ is $R$-trivial.

In our approach to Hilbert's 90 for symbols one needs a
parameterization of the splitting fields of symbols.

{\bf Definition.} \it
  Let $u=\{a_1,\dots,a_n\}\bmod p$ be a symbol.  A \emph{$p$-generic
  splitting variety} for~$u$ is a smooth variety~$X$ over~$k$ such
  that for every splitting field~$F$ of~$u$ there exists a finite
  extension $F'/F$ of degree prime to~$p$ and a morphism $\Spec F'\to
  X$.
\rm

{\bf Theorem.} \it
  Suppose $\car k=0$.  Let $m\geq 3$ and suppose for $n\leq m$ and
  every symbol $u=\{a_1,\dots,a_n\}\bmod p$ over all fields over~$k$
  there exists a $p$-generic splitting variety for~$u$ of dimension
  $p^{n-1}-1$.  Then Hilbert's 90 holds for such symbols.
\rm

The proof of this theorem is outlined in \cite{Rost:98x}.

For $n=2$ one can take here the Severi-Brauer varieties and for
$n=3$ the Merkurjev-Suslin varieties.  Hence we have:

{\bf Corollary.} \it
  Suppose $\car k=0$.  Then Hilbert's 90 holds for symbols of
  weight~$\leq 3$.
\rm

\relax

\bibliographystyle{amsplainR}

\raggedbottom \providecommand{\bibextra}{} \bibextra
  \providecommand{\WWWaddress}[1]{#1}
  \providecommand{\REM}[1]{\marginpar{\tt\raggedright#1}}
  \providecommand{\mhy}{-} \providecommand{\cyr}{} \providecommand{\cprime}{'}
  \providecommand{\SortNoop}[1]{} \providecommand{\transl}[1]{[#1]}
\providecommand{\bysame}{\leavevmode\hbox
to3em{\hrulefill}\thinspace}

\label{lastpage}

\end{document}